\documentclass[draft, 10pt]{article}
\usepackage{amssymb, amsthm, amsmath}
\numberwithin{equation}{section}

\newcommand{\F}{\mathbb{F}}

\newcommand{\Z}{\mathbb{Z}}

\newcommand{\Q}{\mathbb{Q}}

\newcommand{\bfP}{\mathbb{P}}
\newcommand{\A}{\mathbb{A}}

\newcommand{\e}{\emph}
\newcommand{\rom}{\mathrm}
\newcommand{\ov}{\overline}
\newcommand{\ma}{\mathbf}

\newcommand{\ben}{\begin{enumerate}}
\newcommand{\een}{\end{enumerate}}
\newcommand{\eit}{\begin{itemize}}
\newcommand{\beq}{\begin{equation}}
\newcommand{\eeq}{\end{equation}}

\newcommand{\ve}{\varepsilon}

\newcommand{\del}{\delta}

\newcommand{\lab}{\label}

\newtheorem{thm}{Theorem}
\newtheorem{lem}{Lemma}
\newtheorem{pro}{Proposition}
\newtheorem{cor}{Corollary}
\newtheorem{con}{Conjecture}

\theoremstyle{definition}
\newtheorem*{ack}{Acknowledgement}

\newcommand{\hcf}{\rom{h.c.f.}}
\renewcommand{\mod}{\hspace{-1mm}\pmod}
\newcommand{\colt}[2]{\genfrac{}{}{0pt}{1}{#1}{#2}}

\newcommand{\x}{\ma{x}}

\newcommand{\bt}{\ma{t}}

\newcommand{\af}{\mathsf{aff}}

\newcommand{\conj}{Conjecture }

\renewcommand{\leq}{\leqslant}
\renewcommand{\geq}{\geqslant}

\begin{document}

\title{The density of  rational points on\\ non-singular hypersurfaces, I}

\author{T.D. Browning$^1$ and D.R. Heath-Brown$^2$\\
\small{$^{1}$\emph{School of Mathematics,  Bristol University, Bristol BS8 1TW}}\\
\small{$^2$\emph{Mathematical Institute,24--29 St. Giles',Oxford OX1 3LB}}\\
\small{$^1$t.d.browning@bristol.ac.uk, 
$^2$rhb@maths.ox.ac.uk}}

\date{}
\maketitle

\begin{abstract}
For any $n \geq 3$, let 
$F \in \Z[X_0,\ldots,X_n]$ be a form of degree
$d\geq 5$ that defines a non-singular hypersurface $X \subset \bfP^{n}$.
The main result in this paper is a proof of the fact that 
the number $N(F;B)$ of $\mathbb{Q}$-rational points on $X$ which have
height at most $B$ satisfies
$$
N(F;B)=O_{d,\varepsilon,n}(B^{n-1+\varepsilon}),
$$
for any $\varepsilon>0$.
The implied constant in this estimate depends at most upon $d, \varepsilon$ and
$n$.  New estimates are also obtained for the number of representations of a
positive integer as the sum of three $d$th powers, and for the
paucity of integer solutions to equal sums of 
like polynomials. 
\end{abstract}

\section{Introduction}

For any $n \geq 3$, let $F \in \Z[X_0,\ldots,X_n]$ be a form of degree
$d\geq 2$ that produces a non-singular hypersurface $F=0$ in $\bfP^{n}$.
The primary purpose of this paper is to study the distribution of
rational points on this hypersurface.  To this end we define 
the quantity
$$
N(F;B)=\#\{\x \in \Z^{n+1}: F(\x)=0, ~\hcf(x_0,\ldots,x_n)=1, ~|\x| \leq B\},
$$
for any $B\geq 1$, where $|\x|$ denotes the norm $\max_{0\leq i\leq n}|x_i|$.  
A simple heuristic argument leads one to expect an asymptotic formula
of the shape
$$
N(F;B)\sim c_F B^{n+1-d},
$$
as $B \rightarrow \infty$, for some non-negative constant $c_F$.
Birch \cite{birch} has shown that this is indeed the case
when $n$ is sufficiently large compared with $d$, and the conjecture of Manin \cite{man}
predicts that this estimate should be true as soon as $n \geq 2d$,
with $n \geq 4$.  When $d=n=3$ the quantity $N(F;B)$ may be dominated by
the presence of rational lines contained in the surface $F=0$, in which
case it is more interesting to study the distribution of rational
points on the Zariski open subset formed by deleting the lines from the surface.
This paper is motivated by the following basic conjecture.

\begin{con}\lab{hb-ns}  Let $\ve>0$ and suppose that $F \in \Z[X_0,\ldots,X_n]$
is a non-singular form of degree $d \geq 2$.
Then we have
$$
N(F;B)=O_{d,\ve,n}(B^{n-1+\ve}).
$$
\end{con}

Throughout our work the implied constant in any estimate is absolute unless 
explicitly indicated otherwise.  In the case of \conj \ref{hb-ns} for example, the
constant is permitted to depend only upon 
$d$, $\ve$ and $n$. 
\conj \ref{hb-ns} is a special case of a conjecture due
to the second author \cite[\conj 2]{annal}, which predicts that the same
estimate should hold under the weaker assumption that $F$ is absolutely irreducible.
Although \conj \ref{hb-ns} is essentially best possible in the case
$d=2$, it is almost certainly not so for $d \geq 3$ and $n \geq 4$.
In fact \conj \ref{hb-ns} might be considered
as an approximation to the conjecture of Batyrev and Manin
\cite{bat}, which  leads one to expect an estimate of the shape 
$N(F;B)=o_F(B^{n-1})$ for such values of $d$ and $n$.

In this paper we return to the techniques introduced in our earlier
study of rational points of bounded height on geometrically
integral hypersurfaces \cite{bhbs}.  Thus if
$F\in\Z[X_0,\ldots,X_n]$ is an absolutely irreducible form of
degree $d\geq 2$, the main result in this latter work \cite[Corollary 2]{bhbs}
implies that
\beq\lab{sand}
N(F;B) \ll_{d,\ve,n}\left\{
\begin{array}{ll}
B^{n-7/4+5/(3\sqrt{3})+\ve}, & d=3,\\
B^{n-1+\ve} + B^{n-5/3+3/(2\sqrt{d})+\ve}, & d\geq 4.
\end{array}
\right.
\eeq
In particular this estimate  
establishes Conjecture \ref{hb-ns} for $d \geq 6$, since any
non-singular form is automatically absolutely irreducible.
The conjecture had already been established in the case $d=2$,  and for any $d\geq 2$
provided that $n\leq 4$.  This comprises the combined work of 
both the first and second authors \cite{3fold,bhb,annal}.
In all other cases (\ref{sand}) is currently the best available estimate.  It turns out however that
the techniques used to prove (\ref{sand}) for absolutely
irreducible forms can be made to yield
a sharper estimate in the setting of non-singular forms.  
It is important to observe in this context that \conj \ref{hb-ns}
appears to become rapidly more difficult as $d$ decreases.  Our goal
in the present paper has been merely to produce a bound which
establishes \conj \ref{hb-ns} for all $d \geq 5$, in as simple a
manner as possible.  It seems likely that the exponent may be improved
with further work, but that new ideas would be needed in order to
prove the conjecture when $d=3$ or $4$.
The  following is our main result.

\begin{thm}\lab{1}
Let $\ve>0$ and suppose that $F \in \Z[X_0,\ldots,X_n]$ is a non-singular form of degree
$d\geq 4$.  Then we have
$$
N(F;B) \ll_{d,\ve,n} 
B^{n-1+\ve} + B^{n-2+2/\sqrt{d}+1/(d-1)-1/((d-2)\sqrt{d})+\ve}.
$$
In particular Conjecture \ref{hb-ns} holds for $d\geq 5$.
\end{thm}

The estimates in Theorem \ref{1} and (\ref{sand}) coincide
when $d=~4$.  It is interesting to compare Theorem \ref{1} with earlier work of the second author
\cite[Theorem 2]{hb-india}. There the bound 
$N(F;B)=O_{\ve,F}(B^{n-1+\ve})$ 
is established for all non-singular
forms $F$ in at least $10$ variables, provided that $F$ has degree $d \geq 3$.   
As indicated the formulation of \cite[Theorem 2]{hb-india} involves
an implied constant that is allowed to depend upon the
coefficients of $F$.  However a straightforward examination
of the proof of this result reveals that the bound is actually uniform
in non-singular hypersurfaces of fixed degree $d$ and fixed dimension $n-1$.
It will be convenient to state this observation formally here.

\begin{thm}\lab{10}
\conj \ref{hb-ns} holds for $n \geq 9$.
\end{thm}

We are now in a position to combine all of our various results that pertain
to \conj \ref{hb-ns}.  We have already observed that this conjecture holds for all quadrics, and all forms in at most $5$
variables.  On combining Theorem \ref{1} and Theorem
\ref{10} with these facts we therefore obtain the following result.

\begin{cor}
\conj \ref{hb-ns} holds for all values of $d$ and $n$ except possibly
for the eight cases in which $d=3,4$ and $n=5,6,7,8$.
\end{cor}

In a future paper we shall use methods from \cite{bhb} to complete
the proof of \conj \ref{hb-ns} for the eight remaining cases.
The proof of Theorem \ref{1} is based upon the proof of
(\ref{sand}),  and it will be useful to recall the basic idea here. 
For any $\nu \geq 2$, let $f \in \ov{\Q}[T_1,\ldots,T_\nu]$ be a 
polynomial of total degree $\delta$, which produces a non-singular hypersurface $f=0$ in
$\A^{\nu}$.  Here we note that an affine variety is said to be
non-singular if and only if it has  a non-singular projective model, and
we shall henceforth think of a polynomial $f$ as being non-singular if
and only if the corresponding hypersurface $f=0$ is non-singular.  With this in
mind we define the counting function 
\begin{equation}\lab{aff-count}
M(f;B)=\#\{\bt \in \Z^\nu: f(\bt)=0, ~|\bt| \leq B\},
\end{equation}
for any $B \geq 1$.
Then the main trick behind the proof of Theorem \ref{1} is the simple observation that
\beq\lab{trick}
N(F;B) \leq \sum_{|b|\leq B} M(f_b;B),
\eeq
where $f_b=f_b(T_1,\ldots,T_{n})$ denotes the polynomial
$F(b,T_1,\ldots,T_n)$.  Typically $f_b$ will be non-singular and have
degree $d$, so that the problem of proving good upper bounds
for $N(F;B)$ is replaced by the problem of proving good upper bounds
for $M(f;B)$ for non-singular polynomials $f \in
\ov{\Q}[T_1,\ldots,T_\nu]$ of degree $\delta\geq 3$. 
The procedure is then to argue by induction on $\nu$, and it turns out
that the key ingredient is a sharp bound for the inductive
base $\nu=3$.  Thus the main work in the proof of Theorem \ref{1} is 
actually taken up with a demonstration of the following result.

\begin{thm}\lab{2}
Let $\ve>0$ and 
suppose that $f \in \Z[T_1,T_2,T_3]$ is a non-singular polynomial of degree
$\delta\geq 4$.  Then we have
$$
M(f;B)\ll_{\delta,\ve}
B^{1+\ve}+B^{2/\sqrt{\del}+1/(\del-1)-1/((\del-2)\sqrt{\del})+\ve}.
$$
\end{thm}

As highlighted above, the implied constant here is only permitted to depend upon
$\del$ and $\ve$, and not on the individual coefficients of $f$.
It is interesting to compare Theorem \ref{2} with the corresponding
result \cite[Theorem 2]{bhbs}, 
that constitutes the inductive base at the heart of (\ref{sand}).
This latter result applies to the class of all degree $\delta$
polynomials in $\Z[T_1,T_2,T_3]$, whose homogeneous part of degree
$\delta$ is absolutely irreducible.  However the estimate
obtained is weaker than that of Theorem \ref{2}, as soon as $\delta\geq 5$.  In
fact the assumption of non-singularity in Theorem \ref{2} affords
us much better control over the curves of low degree contained in the
surface $f=0$, and this is crucial to the proof.
Theorem \ref{2} will be established in \S \ref{surfaces}, and this
will then be used in \S \ref{hypersurfaces} to deduce Theorem \ref{1}.  

It turns out that Theorem \ref{2} has applications to a rather different
set of problems.  We begin by seeing how it can be used to deduce
better estimates for the number of representations of a given integer
as the sum of three powers.
For given integers $d \geq 2$ and $N \geq 1$, define $r_d(N)$ to be the
number of positive integers $t_1,t_2,t_3$ such that 
\begin{equation}\lab{3powers-eqn}
t_1^d+t_2^d+t_3^d=N.
\end{equation}
It is plain that this equation defines a non-singular affine surface. 
Until recently the best available estimate for $r_d(N)$ was the trivial estimate
\beq\lab{rN-triv}
r_d(N)=O_{d,\ve}(N^{1/d+\ve}),
\eeq
whose exponent is of course
essentially best possible in the case $d=2$.  This was improved for $d
\geq 8$ by the second author \cite[Theorem 13]{annal}, who has replaced the
exponent $1/d$ by $\theta/d$ for $d \geq 2$, with
\beq\lab{exp-hb}
\theta=2/\sqrt{d}+2/(d-1).
\eeq
By combining the proof of this result with the key step in the proof
of Theorem~\ref{2}, we shall establish the following result in \S \ref{cors}.

\begin{cor}\lab{cor1}
Let $\ve>0$ and suppose that $d \geq 4$.  Then we have 
$$
r_d(N) =O_{d,\ve} (N^{\theta/d+\ve}),
$$
with
$$
\theta= 2/\sqrt{d}+1/(d-1)-1/((d-2)\sqrt{d}).
$$
\end{cor}

It is clear that Corollary \ref{cor1} improves upon the trivial
estimate (\ref{rN-triv}) for $d \geq 5$, and upon (\ref{exp-hb}) in
every case. In particular this provides
further evidence for the well-known conjecture that 
$
r_d(N)=O_{d,\ve}(N^{\ve}),
$ 
as soon as $d \geq 4$.

Theorem \ref{2} also has applications to paucity problems for equal
sums of like polynomials.  Let $f \in \Z[T]$ be a polynomial of degree
$d \geq 3$, and let $B$ be a positive integer.  Then for any $s \geq 2$ we define  $L_s(f;B)$ to be the number 
of positive integers $x_1,\ldots,x_{2s} \leq B$ such that 
\beq\lab{wedding}
f(x_1)+\cdots+f(x_s)=f(x_{s+1})+\cdots+f(x_{2s}). 
\eeq
It is conjectured that
$L_s(f;B)$ is dominated by the $s!B^s$ trivial solutions, in which
$x_1,\ldots,x_s$ are a permutation of $x_{s+1},\ldots,x_{2s}$.
Outside of work due to Linnik \cite{linnik},  and later the first author \cite{polys},
there has been relatively little attention paid to the quantity
$L_s(f;B)$ for arbitrary polynomials $f \in \Z[T]$ of degree $d$.  In
this setting, when $s=2$, it follows from \cite[Theorem 1]{polys} that 
$$
L_2(f;B) = 2B^2(1+o_f(1)),
$$ 
provided that $d \geq 7$.  
Earlier work of Wooley \cite{wooley} establishes such an estimate in
the case $d=3$.
We are now ready to present the improvements that we have been able to
make to this problem, the details of which may be found in \S \ref{cors}.

\begin{cor}\lab{4waring}
Let $\ve>0$ and suppose that $f \in \Z[T]$ is a polynomial of degree $d$.  Then we have
$$
L_2(f;B) = 2B^2(1+o_f(1))
$$ 
for $d \geq 5$.
\end{cor}

\begin{ack}
While working on this paper, the first
author was supported at Oxford University by
EPSRC grant number GR/R93155/01.
\end{ack}

\section{Affine surfaces}\lab{surfaces}

In this section we shall prove  Theorem \ref{2}.  
Consequently, let $f \in \Z[T_1,T_2,T_3]$  be a non-singular
polynomial of degree $\del \geq 4$.  
Throughout the proof it will be convenient to work
with the projective model for the affine surface $f=0$. Thus let
$X\subset \bfP^3$ be the non-singular surface defined by the form
\begin{equation}\lab{F}
F(X_0,X_1,X_2,X_3)=X_0^\del f(X_1/X_0,X_2/X_0,X_3/X_0),
\end{equation}
of degree $\del$.  We may clearly assume that $F$ is primitive, so
that the highest common factor of its coefficients is $1$. 
Throughout our work we shall always follow the convention that any point $x\in\bfP^3(\Q)$ is 
represented by a vector $\x=(x_0,x_1,x_2,x_3) \in \Z^4$ such that
$\hcf(x_0,x_1,x_2,x_3)=1$, and we shall write $x=[\x] \in \bfP^3(\Z)$
to express this fact. With this in mind we define the counting function 
$$
N^\af(\Sigma;B)=\#\{[1,x_1,x_2,x_3] \in \Sigma\cap \bfP^{3}(\Z): 
~H([1,x_1,x_2,x_3]) \leq B\},
$$
for any $B\geq 1$ and any locally closed subset $\Sigma\subseteq X$ defined over $\ov{\Q}$.
In particular we clearly have $M(f;B)=N^\af(X;B)$ in the statement of
Theorem \ref{2}.

Let $F\in \Z[X_0,X_1,X_2,X_3]$ be a primitive form of
degree $d \geq 4$ that defines a non-singular surface $X \subset
\bfP^3$, and let $k\geq 2$ be an integer.  Then it will be convenient to introduce the open
subset $X^{(k)} \subseteq X$ which is obtained by deleting all of the
curves from $X$ that are defined over $\overline{\Q}$ and have degree
at most $k$. In particular $X^{(k)}$ is plainly
empty for $k \geq d$ and so we henceforth make the assumption that
\beq\lab{k}
2\leq k \leq d-1.
\eeq  
With this notation in mind our first goal is to establish the
following estimate.

\begin{pro}\lab{pro-suffice}
We have 
$$
N^\af(X^{(k)};B)\ll_{d,\ve}
B^{2/\sqrt{d}+1/(k+1)-1/(k\sqrt{d})+\ve}, 
$$
for any $\ve>0$.
\end{pro}

Before proceeding with this task, we first indicate how this will
suffice to complete the proof of Theorem \ref{2}.
By a result of Colliot-Th\'el\`ene
\cite[Appendix]{annal}, we see that the non-singular surface $X$ 
contains $O_d(1)$ geometrically integral curves
of degree at most $d-2$.   We now need the following estimate due to
Pila \cite[Theorem A]{pila}.

\begin{lem}\lab{pilaA}
Let $\ve>0$ and suppose that $C\subset \A^\nu$ is a geometrically integral curve of degree
$\delta \geq 1$.  Then we have
$$
\#\{\ma{t} \in C\cap \Z^\nu: ~|\bt| \leq B\}
=O_{\delta,\ve,\nu}(B^{1/\delta+\ve}).
$$
\end{lem}

Hence it follows from Lemma \ref{pilaA} that the overall contribution to $N^\af(X;B)$
from the set of curves of degree at most $d-2$ is
$O_{d,\ve}(B^{1+\ve})$, since we have already seen that there are only
$O_d(1)$ such curves contained in $X$. On combining this with an application of
Proposition \ref{pro-suffice} with $k=d-2$ we therefore obtain the estimate
$$
N^\af(X;B)\ll_{d,\ve}
B^{1+\ve}+B^{2/\sqrt{d}+1/(d-1)-1/((d-2)\sqrt{d})+\ve}.
$$
This completes the deduction of Theorem \ref{2} from Proposition \ref{pro-suffice}.

We proceed by establishing Proposition \ref{pro-suffice}.
For any prime $p$ we shall write $X_p$ for the surface defined over
$\F_p$ that is obtained by
reducing the coefficients of $F$ modulo $p$, and we denote 
the set $X_p\cap\bfP^3(\F_p)$ by $X_p(\F_p)$.
It will be convenient to define 
$$
S(\Sigma;B)=\{[1,x_1,x_2,x_3] \in \Sigma\cap \bfP^{3}(\Z): 
~H([1,x_1,x_2,x_3]) \leq B\},
$$ 
for any locally closed subset $\Sigma\subseteq X$ defined over
$\ov{\Q}$, so that in particular $N^\af(\Sigma;B)=\#S(\Sigma;B)$.
Now let $\pi=[1,\pi_1,\pi_2,\pi_3]\in X_p(\F_p)$, where 
$\pi_1,\pi_2,\pi_3$ are always assumed to be in $\F_p$.
We also define the set 
$$
S_p(\Sigma;B,\pi)=\Big\{[1,x_1,x_2,x_3] \in \Sigma\cap \bfP^{3}(\Z): 
\begin{array}{l}
H([1,x_1,x_2,x_3]) \leq B,\\ 
x_i\equiv \pi_i \mod{p}, ~(1\leq i \leq 3)
\end{array}
\Big\},
$$
for any locally closed subset $\Sigma\subseteq X$ defined over $\ov{\Q}$.
We shall base our proof of Proposition \ref{pro-suffice} on the proof of
\cite[Theorem 2]{bhbs}.  
Suppose first that $\log \|F\| \gg_{d}\log B$, with a suitably large
implied constant, where $\|F\|$ denotes the maximum modulus of the
coefficients of $F$.  Then an application of \cite[Lemma 5]{bhbs} reveals
that $N^\af(X^{(k)};B)\leq N^\af(Y;B)$, for some curve $Y\subset
\bfP^3$ of degree $O_d(1)$ that does not contain any geometrically
integral components of degree less than $k$.
The contribution from each such component is clearly 
$O_{d,\ve}(B^{1/(k+1)+\ve})$  by Lemma \ref{pilaA}, and so it suffices
to assume that 
\beq\lab{height-F}
\log \|F\| =O_{d}(\log B)
\eeq
in our proof of  Proposition~\ref{pro-suffice}.
The next step in the argument involves introducing the 
set of $x\in X$ that have multiplicity at most
$2$ on the tangent plane section $X\cap \mathbb{T}_x(X)$.  
We shall henceforth denote this set by $U$, and write $U^{(k)}$ for the
subset of points that do not lie on any curve of degree at most $k$
contained in $X$.    
In particular an application of \cite[Lemmas 10 and 11]{bhbs} 
implies that $U$ is a non-empty open subset of $X$, and once combined
with Lemma \ref{pilaA} we may deduce that
$$
N^\af(X^{(k)};B)=N^\af(U^{(k)};B)+O_{d,\ve}(B^{1/(k+1)+\ve}).
$$
In summary it will therefore suffice to establish 
Proposition~\ref{pro-suffice} under the assumption that
(\ref{height-F}) holds, and with $X^{(k)}$ replaced by $U^{(k)}$.
For any prime $p$ we shall define $U_p$ to be
the open set of non-singular
points on $X_p$ which have multiplicity at most $2$ on the tangent plane
section at the point.
Our main tool in the proof of Theorem \ref{2} is the following
result \cite[Lemma 12]{bhbs}.

\begin{lem}\lab{14}
Let $\ve>0$ and let $U, X$ be as above. Then there exists 
a set $\Pi$ of $O_{d,\ve}(1)$ primes $p$, with
\begin{equation}\lab{prime}
B^{1/{\sqrt{d}}+\ve}\ll_{d,\ve} p\ll_{d,\ve} B^{1/{\sqrt{d}}+\ve},
\end{equation}
such that the following holds.
For each $\pi=[1,\pi_1,\pi_2,\pi_3]\in U_p(\F_p)$, there exists a form
$G_\pi \in \Z[X_0,X_1,X_2,X_3]$  
of degree $O_{d,\ve}(1)$ which is not divisible by $F$, such that
$$
S(U;B)=
\bigcup_{p\in \Pi}\bigcup_{\pi\in U_p(\F_p)}
\{[1,x_1,x_2,x_3] \in
S_p(U;B,\pi): G_\pi(1,x_1,x_2,x_3)=0\}.
$$
\end{lem}

Applying Lemma \ref{14} we therefore deduce that there exists a set $\Pi$ of
$O_{d,\ve}(1)$ primes $p$, with (\ref{prime}) holding, such that
$$
N^\af(U^{(k)};B)
\leq
\sum_{\colt{p\in \Pi}{\pi\in U_p(\F_p)}}
\#\{[1,x_1,x_2,x_3] \in
S_p(U;B,\pi): G_\pi(1,x_1,x_2,x_3)=0\}.
$$
Here the $G_\pi \in \Z[X_0,X_1,X_2,X_3]$  are a finite set of forms
indexed by points $\pi=[1,\pi_1,\pi_2,\pi_3]\in U_p(\F_p)$. 
For each $\pi\in U_p(\F_p)$, the form 
$G_\pi$ has degree $O_{d,\ve}(1)$ and is not divisible by $F$.  
For any integral component $Y=Y_\pi \subset \bfP^3$ of the curve
$F=G_\pi=0$ we write 
$$
N_p^\af(Y;B,\pi)=\#S_p(Y;B,\pi).
$$
Now let $Y_1, \ldots, Y_s$ be the collection of those geometrically integral components of
$F=G_\pi=0$ which either have dimension $0$, or have dimension $1$ and
degree at least $k+1$. Then we will have
$s=O_{d,\ve}(1)$ and
$$
\#\{[1,x_1,x_2,x_3] \in
S_p(U^{(k)};B,\pi): G_\pi(1,x_1,x_2,x_3)=0\}
\leq \sum_{1\leq j\leq s} N_p^\af(Y_j;B,\pi).
$$

As we vary over primes $p\in \Pi$ and points 
$\pi\in U_p(\F_p)$, let $I$ denote the set of components 
of the curves $F=G_\pi=0$ that have dimension $0$.  
Since $\#U_p(\F_p) \leq \#X_p(\F_p)=O_d(p^2)$ for each $p\in \Pi$, it
easily follows from (\ref{prime}) that
$$
\#I \leq s\#\Pi. \max_p\#U_p(\F_p) =O_{d,\ve}(B^{2/\sqrt{d}+2\ve}).
$$
Thus the overall contribution to $N^\af(U^{(k)};B)$ from the set $I$ is
satisfactory for Proposition \ref{pro-suffice}.
We henceforth fix a choice of prime $p\in \Pi$ and a point 
$\pi\in U_p(\F_p)$. Let $G_\pi$ be the corresponding form produced by
Lemma \ref{14} and let 
$Y=Y_\pi \subset \bfP^3$ be any geometrically integral component of the curve
$F=G_\pi=0$ of dimension $1$ and degree $e\geq k+1$. 
We now employ the following estimate \cite[Proposition~2]{bhbs}.

\begin{lem}\lab{pro}
Assume that $Y$ has degree $e\geq 3$.  Then we have 
$$
N_p^\af(Y;B,\pi)\ll_{d,e,\ve} B^{1/e-1/((e-1)\sqrt{d})}.
$$
\end{lem}

As above, we note that $\#U_p(\F_p)=O_{d,\ve}(B^{2/\sqrt{d}+2\ve})$ for each $p\in \Pi$.
On observing that 
$$
3\leq k+1 \leq e =O_{d,\ve}(1),
$$  
by (\ref{k}), we may therefore take $e\geq k+1$ in
Lemma \ref{pro}  in order to conclude the proof that
$$
N^\af(U^{(k)};B)\ll_{d,\ve} B^{2/\sqrt{d}+1/(k+1)-1/(k\sqrt{d})+2\ve}.
$$
This suffices for the proof of Proposition \ref{pro-suffice}, and so
also for that of Theorem~\ref{2}.

\section{Affine hypersurfaces}\lab{hypersurfaces}

In this section we shall establish Theorem \ref{1}.    Before doing so
we take a moment to record a ``trivial'' upper bound for the quantity
$N(G;B)$, which has the advantage of applying under the
sole assumption that the form $G$ is non-zero.
The following result is due to the second author \cite[Theorem 1]{annal}.

\begin{lem}\lab{triv}
Let $G \in \Z[X_0,\ldots,X_n]$ be a non-zero form of
degree $d$.  Then we have
$$
N(G;B) =O_{d,n}(B^{n}).
$$
\end{lem}

As outlined in the introduction, the main ingredient in
the proof of Theorem~\ref{1} will be Theorem \ref{2}.
Let $n \geq 3$ and let $F \in \Z[X_0,\ldots,X_n]$ be a non-singular form of
degree $d\geq 4$.   Then any affine
model of $F$ is also non-singular and has degree $d$.  
For any $b \in \Z$ we define the polynomial
$$
f_b(X_1,\ldots,X_{n})=F(b,X_1,\ldots,X_n).
$$
By Lemma \ref{triv} we clearly have
$N(F(0,X_1,\ldots,X_n);B)=O_{d,n}(B^{n-1}).$  Hence, 
in view of (\ref{trick}), we deduce that
\beq\lab{step1-1}
N(F;B) \ll_{d,n} B^{n-1}+ \sum_{0<|b|\leq B} M(f_b;B),
\eeq
where $f_b \in \Z[X_1,\ldots,X_{n}]$ is a non-singular polynomial of degree
$d$ for $b \neq 0$.
The primary aim of this section is to establish the following result,
which may be of independent interest.

\begin{pro}\lab{pilaA-1}
Let $\ve>0$ and let $\nu\geq 3$.  Suppose that $f \in \Z[T_1,\ldots,T_\nu]$ is a non-singular polynomial of degree
$\del\geq 4$.  Then we have
$$
M(f;B) \ll_{\del,\ve,\nu} B^{\nu-2+\ve}+B^{\nu-3+2/\sqrt{\del}+1/(\del-1)-1/((\del-2)\sqrt{\del})+\ve}.
$$
\end{pro}

Before proceeding with the proof of Proposition \ref{pilaA-1}, we
first see how it suffices to complete the proof of Theorem
\ref{1}.   But this is simply a matter of using (\ref{step1-1}) to
deduce that
\begin{align*}
N(F;B) &\ll_{d,\ve,n} B^{n-1}+ \sum_{0<|b|\leq B} \big(B^{n-2+\ve} +
B^{n-3+2/\sqrt{d}+1/(d-1)-1/((d-2)\sqrt{d})+\ve}\big)\\ 
&\ll_{d,\ve,n} B^{n-1+\ve} + B^{n-2+2/\sqrt{d}+1/(d-1)-1/((d-2)\sqrt{d})+\ve},
\end{align*}
for any $\ve>0$. 
This completes the deduction of Theorem \ref{1} from Proposition
\ref{pilaA-1}.
The proof of Proposition \ref{pilaA-1} will involve taking repeated
hyperplane sections of the non-singular affine hypersurface $f=0$.  
We begin by establishing the following auxiliary result, which will
help us to control certain bad hyperplane sections that arise in our argument.

\begin{lem}\lab{elim-1}
Let  $f \in \Z[T_1,\ldots,T_\nu]$ be a non-singular polynomial of degree $\del$.  Then there
exists  a matrix $A=\{a_{ij}\}_{i,j\leq \nu} \in SL_\nu(\Z)$ with
$\max_{i,j}|a_{ij}| \ll_{\del,\nu} 1$, such that 
if $g(S_1,\ldots,S_\nu)=f(A^{-1}\ma{S})$ then the following holds:
\ben 
\item
$g(S_1,\ldots,S_\nu)$ is non-singular and has degree $\del$.
\item
There exists $k \in \Q$ such that $g(k,S_2,\ldots, S_\nu)$ is
non-singular and has degree $\del$.
\een
\end{lem}

\begin{proof}
The first part of Lemma \ref{elim-1} is
true for any polynomial obtained by non-singular linear
transformation from $f$. 
To establish the second part of the
lemma, we consider the degree $\del$ projective model $F(\ma{X})=F(X_0,\ldots,X_{\nu})$ for
$f$.  Then since $F$ is non-singular, there exists a
non-zero form $\hat{F}(\ma{Y})=\hat{F}(Y_0,\ldots,Y_{\nu})$ of degree $O_{\del,\nu}(1)$, such that
$\hat{F}(\ma{y})=0$ whenever the pair of equations
$$
F(\ma{X})=\ma{X}.\ma{y}=0
$$
produce a singular form.  
In fact the theory of the
dual variety ensures that $\hat{F}$ is
irreducible and has degree at least $ 2$, although we shall not
need these facts here.  
It follows from Lemma \ref{triv} that we may find a vector $\ma{b} \in
\Z^{\nu+1}$ such that  $|\ma{b}| \ll_{\del,\nu} 1$ and 
$b_1\hat{F}(\ma{b}) \neq 0.$  Set $h=\hcf(b_1,\ldots,b_\nu)$ and write
$$
k=-h^{-1}b_0, \quad a_{1j}=h^{-1}b_j,
$$ 
for $1 \leq j \leq \nu$.  Then 
it follows that $\hcf(a_{11},\ldots, a_{1\nu})=1$, and 
we may therefore find a matrix $A=\{a_{ij}\}_{i,j\leq \nu} \in GL_\nu(\Z)$ such that $\det(A)=1$ and $\max_{i,j}|a_{ij}|
\ll_{\del,\nu} 1$.  We then see that the variety
$g(k,S_2,\ldots,S_\nu)=0$ is alternatively given by
$f(A^{-1}(k,S_2,\ldots,S_\nu))=0$, and so is equal to the variety
$$
f(T_1,\ldots,T_\nu)=0, \quad a_{11}T_1+\cdots +a_{1\nu} T_\nu =k,
$$
which is non-singular and has degree $\delta$ by construction.  This completes the
proof of Lemma~\ref{elim-1}.
\end{proof}

We now proceed with the proof of Proposition \ref{pilaA-1}, for which
we employ induction on $\nu$.  The inductive base $\nu=3$ is
satisfactory by Theorem \ref{2}.
By Lemma \ref{elim-1} there exists a non-singular polynomial
$g=g(S_1,\ldots,S_\nu)$  and a constant $c \ll_{\del,\nu} 1$ such that
$$
M(f;B)\leq M(g;c B)= \sum_{|\kappa| \leq cB} M(g_\kappa; cB),
$$
where we write 
$g_\kappa(S_2,\ldots,S_{\nu})=g(\kappa,S_2,\ldots,S_\nu)$ for any fixed choice of $\kappa$.
Then if $g_\kappa$ is non-singular and has degree $\del$, we may employ the
induction hypothesis to deduce that 
$$
M(g_\kappa;c B)\ll_{\del,\ve,\nu} B^{\nu-3+\ve}+B^{\nu-4+2/\sqrt{\del}+1/(\del-1)-1/((\del-2)\sqrt{\del})+\ve}.
$$
Alternatively, if $g_\kappa$ is singular or has degree is less than
$\delta$, we use the trivial estimate $M(g_\kappa; c
B)=O_{\del,\nu}(B^{\nu-2})$ coming from \cite[Equation (2.3)]{bhb}.  Note that
$g_\kappa$ cannot vanish identically since $g$ is an irreducible polynomial
of degree at least two.  We
claim that that there are just $O_{\del,\nu}(1)$ values of $\kappa \in \Z$ for which
$g_\kappa$ is a bad hyperplane section, where we think of a hyperplane
section as being bad if it produces a polynomial which is either
singular or of strictly smaller degree.  To see the claim we note 
as in the proof of Lemma \ref{elim-1} that there is a polynomial
$\hat{g}(\kappa)$ of degree $O_{\del,\nu}(1)$, which vanishes precisely when $g_\kappa$
is a bad hyperplane section. But Lemma \ref{elim-1} also ensures that
$\hat{g}$ does not vanish identically, which therefore establishes the
claim.  Hence
\begin{align*}
M(f;B)&\ll_{\del,\ve,\nu} B^{\nu-2}+ \sum_{|\kappa| \leq cB}
\big(B^{\nu-3+\ve}+B^{\nu-4+2/\sqrt{\del}+1/(\del-1)-1/((\del-2)\sqrt{\del})+\ve}\big) \\
&\ll_{\del,\ve,\nu} B^{\nu-2+\ve}+B^{\nu-3+2/\sqrt{\del}+1/(\del-1)-1/((\del-2)\sqrt{\del})+\ve}.
\end{align*}
This establishes Proposition  \ref{pilaA-1}, and so
completes the proof of Theorem~\ref{1}.

\section{Proof of Corollaries \ref{cor1} and \ref{4waring}}\lab{cors}

We begin by considering the affine surface (\ref{3powers-eqn}), with
a view to establishing Corollary \ref{cor1}.
Thus it plainly suffices to show that 
for any $B \geq N^{1/d}$ we have 
\begin{equation}\lab{goal-3powers}
M(h;B)\ll_{d,\ve} B^{2/\sqrt{d}+1/(d-1)-1/((d-2)\sqrt{d})+\ve},
\end{equation}
in the notation of (\ref{aff-count}), where $h$ denotes the non-singular polynomial
$$
h(T_1,T_2,T_3)=T_{1}^{d}+T_{2}^{d}+T_{3}^{d}-N.
$$
Let $S\subset \A^3$ denote the surface $h=0$.  Then a straightforward inspection of
the proof of \cite[Theorem 13]{annal} shows that the overall contribution to $M(h;B)$ from the
union of curves of degree at most $d-2$ contained in $S$, is 
$O_{d, \ve}(B^{1/d+\ve})$.  This is plainly satisfactory for (\ref{goal-3powers}).
In order to complete the proof of Corollary \ref{cor1} it
therefore suffices to apply Proposition \ref{pro-suffice} with $k=d-2$, much as
in the proof of Theorem \ref{2}.

Finally we establish Corollary \ref{4waring}.  For this we shall take
$k=d-2$ in the statement of Proposition \ref{pro-suffice}.  But then
it follows that Hypothesis $[d,\theta_d]$ holds
in \cite[Theorem 3]{polys}, with
$\theta_d=2/\sqrt{d}+1/(d-1)-1/((d-2)\sqrt{d})$. Let $B$ be a positive
integer and let $L_s^{(0)}(f;B)$ denote the contribution to $L_s(f;B)$ from
the integer solutions to (\ref{wedding}) in which 
$x_1,\ldots,x_s$ are not a permutation of $x_{s+1},\ldots,x_{2s}$.
Then we deduce that 
$$
L_s(f;B)=s!B^s+ L_s^{(0)}(f;B),
$$ 
with 
$$
L_s^{(0)}(f;B) \ll_{\ve,f} B^{2s-3+\ve}\Big(B^{1/3} +
  B^{2/\sqrt{d}+1/(d-1)-1/((d-2)\sqrt{d})} \Big).
$$
This therefore establishes Corollary \ref{4waring}.

\end{document}